\documentclass{amsart}
\usepackage[latin1]{inputenc}

\usepackage{amssymb,amsmath}
\usepackage{amsthm}
\usepackage{paralist}
\usepackage{url}

\theoremstyle{plain}
\newtheorem{theorem}{Theorem}
\newtheorem{proposition}[theorem]{Proposition}
\newtheorem{lemma}[theorem]{Lemma}

\theoremstyle{remark}
\newtheorem{remark}{Remark}
\newtheorem{claim}{Claim}

\newcommand{\vect}[1]{\ensuremath{\mathbf{#1}}} 
\newcommand{\card}[1]{\ensuremath{\lvert{#1}\rvert}} 
\newcommand{\minor}[3]{\ensuremath{{#1}_{{#2} \gets {#3}}}} 
\DeclareMathOperator{\ess}{ess} 
\DeclareMathOperator{\essl}{ess^{<}} 
\DeclareMathOperator{\gap}{gap} 
\DeclareMathOperator{\qa}{qa} 
\DeclareMathOperator{\id}{id} 

\newcommand{\oddsupp}{\ensuremath{\mathrm{oddsupp}}} 

\newcommand{\Aneq}[1][n]{\ensuremath{A^{#1}_{=}}}

\begin{document}

\title[Generalizations of \'{S}wierczkowski's lemma and the arity gap]{Generalizations of \'{S}wierczkowski's lemma and the arity gap of finite functions}

\author{Miguel Couceiro}
\address[M. Couceiro]{Department of Mathematics \\
University of Luxembourg \\
162a, avenue de la Faïencerie \\
L--1511 Luxembourg}
\email{miguel.couceiro@uni.lu}
\author{Erkko Lehtonen}
\address[E. Lehtonen]{Department of Combinatorics and Optimization \\
University of Waterloo \\
200 University Avenue West \\
Waterloo, Ontario, N2L~3G1 \\
Canada} 
\address{Department of Mathematics \\
Tampere University of Technology \\
P.O. Box 553 \\
FI--33101 Tampere \\
Finland}
\email{erkko.lehtonen@uni.lu}
\thanks{This work was partially supported by the Academy of Finland, grant \#120307}
\begin{abstract}
\'{S}wierczkowski's lemma---as it is usually formulated---asserts that if $f \colon A^n \to A$ is an operation on a finite set $A$, $n \geq 4$, and every operation obtained from $f$ by identifying a pair of variables is a projection, then $f$ is a semiprojection. We generalize this lemma in various ways. First, it is extended to $B$-valued functions on $A$ instead of operations on $A$ and to essentially at most unary functions instead of projections. Then we characterize the arity gap of functions of small arities in terms of quasi-arity, which in turn provides a further generalization of \'{S}wierczkowski's lemma. Moreover, we explicitly classify all pseudo-Boolean functions according to their arity gap. Finally, we present a general characterization of the arity gaps of $B$-valued functions on arbitrary finite sets $A$.
\end{abstract}
\maketitle

\section{Introduction}

\'{S}wierczkowski's lemma has fundamental consequences in universal algebra (see, e.g., \cite{FKM,PR,PRS,Rosenberg}). It is usually formulated as follows: given an operation $f \colon A^n \to A$, $n \geq 4$, if every operation obtained from $f$ by identifying a pair of variables is a projection, then $f$ is a semiprojection, i.e., there is a $t \in \{1, \dotsc, n\}$ such that $f(a_1, \dotsc, a_n) = a_t$ whenever $a_i = a_j$ for some $1 \leq i < j \leq n$. The importance of this result in clone theory was made apparent in several studies, in particular in Rosenberg's classification of minimal clones \cite{Rosenberg}. This classification was extended to clones of multioperations in \cite{PR}, where \'{S}wierczkowski's lemma was accordingly adjusted for multiprojections in the multi-valued case.

\'{S}wierczkowski's lemma can be generalized in terms of a quasi-ordering of functions $f \colon A^n \to B$, the so-called simple minor relation \cite{CL,CP}: a function $f$ is a simple minor of a function $g$ if $f$ can be obtained from $g$ by permutation of variables, addition or deletion of inessential variables, and identification of variables (see Section \ref{sec:minors}). With this terminology, in Section \ref{sec:quasi} we extend \'{S}wierczkowski's lemma to the following: if every function obtained from $f \colon A^n \to B$ by identifying a pair of variables is essentially unary and $n = 2$ or $n \geq 4$ (constant and $n \geq 2$, respectively), then there exists an essentially unary function (a constant function, respectively) $g \colon A^n \to B$ such that $f|_{\Aneq} = g|_{\Aneq}$, where 
\[
\Aneq = \{(a_1, \dotsc, a_n) \in A^n : \text{$a_i = a_j$ for some $1 \leq i < j \leq n$}\}.
\]

The simple minor relation plays an important role in the equational approach to function class definability. Ekin \emph{et al.} \cite{EFHH} showed that the equationally definable classes of Boolean functions coincide with the initial segments of this quasi-ordering. This result was extended to functions defined on arbitrary, possibly infinite domains in \cite{CF}. These results motivated a study \cite{CP} of the simple minor relation where the structure of the corresponding quasi-ordered set was investigated and where the notion of arity gap had useful applications.

The arity gap of a function $f \colon A^n \to B$ which depends on all of its variables, $n \geq 2$, is defined as the minimum decrease in the number of essential variables when variables of $f$ are identified. The arity gap of a function is obviously at least $1$, and it can be as large as $\card{A}$, as shown by the following example due to Salomaa \cite{Salomaa}. Assume that $\card{A} = k \geq 2$, let $\vect{a} = (a_1, \dotsc, a_k) \in A^k$ be a fixed $k$-tuple such that $a_i \neq a_j$ whenever $i \neq j$, and let $b$ and $c$ be distinct elements of $A$. Define the operation $f \colon A^k \to A$ as follows:
\[
f(\vect{x}) =
\begin{cases}
b, & \text{if $\vect{x} = \vect{a}$,} \\
c, & \text{otherwise.}
\end{cases}
\]
All $k$ variables are essential in $f$, and whenever any two variables are identified, the resulting function is constant, having no essential variables. Thus the arity gap of $f$ is $k$.

Salomaa \cite{Salomaa} showed that the arity gap of any Boolean function is at most $2$. This result was extended to functions defined on arbitrary finite domains by Willard \cite{Willard}, who showed that the same upper bound holds for the arity gap of any function $f \colon A^n \to B$, provided that $n > \card{A}$. In fact, he showed that if the arity gap of such a function $f$ is $2$, then $f$ is totally symmetric. In Section \ref{sec:gap} we consider the arity gap of functions $f \colon A^n \to B$ of small arities, namely, where $n \leq \card{A}$. We give a characterization of such functions by introducing the notion of quasi-arity, which in turn leads to a further generalization of \'{S}wierczkowski's lemma.

In \cite{CL} we strengthened Salomaa's \cite{Salomaa} result on the upper bound for the arity gap of Boolean functions by completely classifying the Boolean functions according to their arity gap. By making use of tools provided by Berman and Kisielewicz \cite{BK} and Willard \cite{Willard}, we obtain in Section \ref{sec:pseudoBoolean} a similar explicit classification of all pseudo-Boolean functions, i.e., functions $f \colon \{0,1\}^n \to B$, where $B$ is an arbitrary set. In Section \ref{sec:concluding}, we present a general characterization of finite functions according to their arity gap which is given in terms of quasi-arity.

\section{Simple variable substitutions}
\label{sec:minors}

Throughout this paper, let $A$ be an arbitrary finite set with $\card{A} = k \geq 2$ elements, and let $B$ be an arbitrary set with at least two elements. A \emph{$B$-valued function (of several variables) on $A$} is a mapping $f \colon A^n \to B$ for some positive integer $n$, called the \emph{arity} of $f$. $A$-valued functions on $A$ are called \emph{operations on $A$.} Operations on $\{0,1\}$ are called \emph{Boolean functions.} For an arbitrary $B$, we refer to $B$-valued functions on $\{0, 1\}$ as \emph{pseudo-Boolean functions.}

For each positive integer $n$, the $n$-ary \emph{projections} $(a_1, \dotsc, a_n) \mapsto a_i$, $1 \leq i \leq n$, are also called \emph{variables} and denoted by $x_i^{(n)}$, or simply by $x_i$ when the arity is clear from the context.
We say that the $i$-th variable $x_i$ is \emph{essential} in $f$, or $f$ \emph{depends} on $x_i$, if there are elements $a_1, \dotsc, a_n, b \in A$ such that
\[
f(a_1, \dotsc, a_{i-1}, a_i, a_{i+1}, \dotsc, a_n) \neq f(a_1, \dotsc, a_{i-1}, b, a_{i+1}, \dotsc, a_n).
\]
In this case, the pair $((a_1, \dotsc, a_{i-1}, a_i, a_{i+1}, \dotsc, a_n), (a_1, \dotsc, a_{i-1}, b, a_{i+1}, \dotsc, a_n))$ is called a \emph{witness of essentiality} of $x_i$ in $f$. If a variable is not essential in $f$, then we say that it is \emph{inessential} in $f$. The number of essential variables in $f$ is called the \emph{essential arity} of $f$, and it is denoted by $\ess f$. If $\ess f = m$, we say that $f$ is \emph{essentially $m$-ary.} Thus the only essentially nullary functions are the constant functions.

We extend the notion of essential variable to partial functions $f \colon S \to B$, where $S \subseteq A^n$. The definition is in fact the same as the one for total functions $A^n \to B$, but now the witnesses of essentiality must be in $S^2$. In other words, we say that the $i$-th variable $x_i$ is \emph{essential} in $f \colon S \to B$, where $S \subseteq A^n$, or $f$ \emph{depends} on $x_i$, if there is a pair
\[
((a_1, \dotsc, a_{i-1}, a_i, a_{i+1}, \dotsc, a_n), (a_1, \dotsc, a_{i-1}, b, a_{i+1}, \dotsc, a_n)) \in S^2,
\]
called a \emph{witness of essentiality} of $x_i$ in $f$, such that
\[
f(a_1, \dotsc, a_{i-1}, a_i, a_{i+1}, \dotsc, a_n) \neq f(a_1, \dotsc, a_{i-1}, b, a_{i+1}, \dotsc, a_n).
\]
With no risk of ambiguity, the notion of essential arity is defined for partial functions as in the case of total functions. 

Suppose that $f \colon S \to B$ is a partial function with $S \subseteq A^n$, and $\ess f = m$. Let us suppose without loss of generality (renaming the variables if necessary) that the essential variables are $x_1, \dotsc, x_m$. Consider the statement
\begin{equation}
\label{st1}
\exists h \colon A^m \to B : \forall (a_1, \dotsc, a_n) \in S : f(a_1, \dotsc, a_n) = h(a_1, \dotsc, a_m).
\end{equation}

Statement \eqref{st1} obviously holds for total functions, i.e., when $S = A^n$. It is not true for partial functions in general. (For a counterexample, assume that $\card{A} \geq 2$, let $S = \{(a, a, \dotsc, a) : a \in A\} \subseteq A^n$, and define $f \colon S \to A$ by $f(a, a, \dots, a) = a$ for all $(a, a, \dotsc, a) \in S$. All variables of $f$ are inessential, yet $f$ is not a constant function.) However, \eqref{st1} is true for partial functions whose domain has a special shape, as described below.

For $n \geq 2$, define the set
\[
\Aneq = \{(a_1, \dotsc, a_n) \in A^n : \text{$a_i = a_j$ for some $i \neq j$}\}.
\]
We define $\Aneq[1] = A$. Note that if $n > \lvert A \rvert$, then $\Aneq = A^n$.

\begin{lemma}
\label{existsh}
Let $f \colon \Aneq \to B$ be a partial function, $n \neq 2$, and $\ess f = m$. Suppose that the essential variables of $f$ are $x_1, \dotsc, x_m$. Then there exists a total function $h \colon A^m \to B$ such that for all $(a_1, \dotsc, a_n) \in \Aneq$, $f(a_1, \dotsc, a_n) = h(a_1, \dots, a_m)$.
\end{lemma}
\begin{proof}
The case $n = 1$ being trivial, we may assume that $n \geq 3$. If $m = n$, then every extension of $f$ into a total function satisfies the required condition. Thus, we can assume that $m < n$. Let $h(a_1, \dotsc, a_m) = f(a_1, \dotsc, a_m, a_m, \dotsc, a_m)$ for all $(a_1, \dotsc, a_n) \in A^n$. Then for any $(a_1, \dotsc, a_n) \in \Aneq$ we have
\begin{align*}
&f(a_1, \dotsc, a_m, a_{m+1}, a_{m+2}, a_{m+3}, \dotsc, a_{n-1}, a_n) && \\
&= f(a_1, \dotsc, a_m, a_m, a_{m+2}, a_{m+3}, \dotsc, a_{n-1}, a_n) && \text{(since $x_{m+1}$ is inessential in $f$)} \\
&= f(a_1, \dotsc, a_m, a_m, a_m, a_{m+3}, \dotsc, a_{n-1}, a_n) && \text{(since $x_{m+2}$ is inessential in $f$)} \\
& \;\, \vdots \\
&= f(a_1, \dotsc, a_m, a_m, a_m, a_m, \dotsc, a_m, a_m, a_n) && \\
&= f(a_1, \dotsc, a_m, a_m, a_m, a_m, \dotsc, a_m, a_m, a_m) && \text{(since $x_n$ is inessential in $f$)} \\
&= h(a_1, \dotsc, a_m). &&
\end{align*}
(Note that the $n$-tuples we used all belong to $\Aneq$.)
\end{proof}

The \emph{composition} of $f \colon B^n \to C$ with $g_1, \dotsc, g_n \colon A^m \to B$ is the function $f(g_1, \dotsc, g_n) \colon A^m \to C$ defined by
\[
f(g_1, \dotsc, g_n)(\vect{a}) = f(g_1(\vect{a}), \dotsc, g_n(\vect{a}))
\]
for all $\vect{a} \in A^m$.

We say that a function $f \colon A^n \to B$ is obtained from $g \colon A^m \to B$ by \emph{simple variable substitution,} or $f$ is a \emph{simple minor} of $g$, if there is a mapping $\sigma \colon \{1, \dotsc, m\} \to \{1, \dotsc, n\}$ such that
\[
f = g(x_{\sigma(1)}^{(n)}, \dotsc, x_{\sigma(m)}^{(n)}).
\]
If $\sigma$ is not injective, then we speak of \emph{identification of variables.} If $\sigma$ is not surjective, then we speak of \emph{addition of inessential variables.} If $\sigma$ is a bijection, then we speak of \emph{permutation of variables.} Observe that each $x_{\sigma(i)}^{(n)}$ is simply an $n$-ary projection, and thus we have that $f$ is a simple minor of $g$ if and only if
\begin{multline*}
\{f(\pi_1, \dotsc, \pi_n) : \text{$\pi_1, \dotsc, \pi_n$ are projections of the same arity}\} \subseteq \\
\{g(\rho_1, \dotsc, \rho_m) : \text{$\rho_1, \dotsc, \rho_m$ are projections of the same arity}\}.
\end{multline*}

From this observation it follows that the simple minor relation constitutes a quasi-order $\leq$ on the set of all $B$-valued functions of several variables on $A$ which is given by the following rule: $f \leq g$ if and only if $f$ is obtained from $g$ by simple variable substitution. If $f \leq g$ and $g \leq f$, we say that $f$ and $g$ are \emph{equivalent,} denoted $f \equiv g$. If $f \leq g$ but $g \not\leq f$, we denote $f < g$. It can be easily observed that if $f \leq g$ then $\ess f \leq \ess g$, with equality if and only if $f \equiv g$. For background, extensions and variants of the simple minor relation, see, e.g., \cite{Couceiro,CP,FH,Lehtonen,LS,Pippenger,Wang,Zverovich}.

\section{Quasi-arity and a generalization of \'{S}wierczkowski's lemma}
\label{sec:quasi}

In this section, we extend \'{S}wierczkowski's lemma to $B$-valued functions on $A$. To this extent, we need to introduce some terminology.

Let $f \colon A^n \to B$, where $n \geq 2$. For indices $i, j \in \{1, \dotsc, n\}$, $i \neq j$, the function $\minor{f}{i}{j} \colon A^n \to B$ obtained from $f \colon A^n \to B$ by the simple variable substitution
\[
\minor{f}{i}{j} = f(x_1^{(n)}, \dotsc, x_{i-1}^{(n)}, x_j^{(n)}, x_{i+1}^{(n)}, \dotsc, x_n^{(n)})
\]
is called a \emph{variable identification minor} of $f$, obtained by identifying $x_i$ with $x_j$.

The \emph{diagonal function} of $f \colon A^n \to B$ is the mapping $\Delta_f \colon A \to B$ defined by $\Delta_f(a) = f(a, a, \dotsc, a)$ for all $a \in A$. Equivalently, in terms of functional composition, $\Delta_f = f(x_1^{(1)}, \dotsc, x_1^{(1)})$.

The proofs of the following two lemmas are straightforward and are left to the reader.

\begin{lemma}
\label{lem:diag}
$\Delta_f = \Delta_{\minor{f}{i}{j}}$ for all $i \neq j$.
\end{lemma}

\begin{lemma}
\label{lem:un}
A function $f \colon A^n \to B$ is essentially at most unary if and only if $f = \Delta_f(x_i^{(n)})$ for some $1 \leq i \leq n$.
\end{lemma}

Let $f \colon A^n \to B$. Any function $g \colon A^n \to B$ satisfying $f|_{\Aneq} = g|_{\Aneq}$ is called a \emph{support} of $f$. The \emph{quasi-arity} of $f$, denoted $\qa f$, is defined as the minimum of the essential arities of the supports of $f$, i.e., $\qa f = \min \ess g$, where $g$ ranges over the set of all supports of $f$. If $\qa f = m$, we say that $f$ is \emph{quasi-$m$-ary.} We call $f$ a \emph{semiprojection,} if there exists a projection that is a support of $f$, in other words, if there is a $t \in \{1, \dotsc, n\}$ such that $f(a_1, \dots, a_n) = a_t$ whenever $a_i = a_j$ for some $1 \leq i < j \leq n$. Note that according to our definition, all projections are semiprojections. Some authors do not consider projections as semiprojections (e.g., \cite{Cs,Quack}).

\begin{remark}
\label{rem1}
\begin{asparaenum}[\rm (i)]
\item
\label{rem:quasi}
If $n > \card{A}$, then quasi-$m$-ary means the same as essentially $m$-ary, in particular, quasi-nullary and quasi-unary mean the same as constant and essentially unary, respectively.

\item
\label{rem:binary}
The quasi-arity of $f \colon A^2 \to B$ is either $0$ or $1$, depending on whether the diagonal function $\Delta_f$ is constant or nonconstant, respectively. Furthermore, the two possible variable identification minors $\minor{f}{1}{2}$ and $\minor{f}{2}{1}$ are equivalent to $\Delta_f$.

\item
\label{rem:unique}
Among all supports of a quasi-nullary function, there is exactly one which is constant. If $f \colon A^2 \to B$ is quasi-unary, then it has exactly two essentially unary supports. For $n \geq 3$, if $f \colon A^n \to B$ is quasi-unary, then among all supports of $f$ there is exactly one which is essentially unary.
\end{asparaenum}
\end{remark}

The following lemma establishes the connection between the quasi-arity of $f$ and the essential arity of the restriction $f|_{\Aneq}$.

\begin{lemma}
\label{qafessf}
For every function $f \colon A^n \to B$, $n \neq 2$, we have $\qa f = \ess f|_{\Aneq}$.
\end{lemma}
\begin{proof}
The case when $n = 1$ being trivial, we may assume that $n \geq 3$. Since any witness of essentiality of an essential variable $x_i$ in $f|_{\Aneq}$ is also a witness of essentiality of $x_i$ in any support of $f$, every essential variable of $f|_{\Aneq}$ must be essential in every support of $f$, and therefore we have $\qa f \geq \ess f|_{\Aneq}$.

Let $\ess f|_{\Aneq} = m$, and assume without loss of generality that the essential variables of $f|_{\Aneq}$ are $x_1, \dotsc, x_m$. By Lemma \ref{existsh}, there is a function $h \colon A^m \to B$ such that for all $(a_1, \dotsc, a_n) \in \Aneq$, $f|_{\Aneq}(a_1, \dotsc, a_n) = h(a_1, \dotsc, a_m)$. By introducing $n - m$ inessential variables, we obtain the function $h' \colon A^n \to B$ given by $h'(a_1, \dotsc, a_n) = h(a_1, \dotsc, a_m)$ for all $(a_1, \dotsc, a_n) \in A^n$, and it is clear that $\ess h' = \ess h$. Then $h'$ is a support of $f$, and since $h$ has arity $m$, we have $\qa f \leq \ess h' = \ess h \leq m$.
\end{proof}

Note that Lemma \ref{qafessf} does not hold for $n = 2$. By Remark \ref{rem1}\eqref{rem:binary}, the quasi-arity of a binary fuction $f \colon A^2 \to B$ is either $0$ or $1$, but $\ess f|_{\Aneq[2]} = 0$.

\begin{lemma}
If a quasi-$m$-ary function $f \colon A^n \to B$ has an inessential variable, then $f$ is essentially $m$-ary.
\end{lemma}
\begin{proof}
The statement is easily seen to hold when $n < 3$, so we may assume that $n \geq 3$. Let $\qa f = m$. By Lemma \ref{qafessf}, $\ess f|_{\Aneq} = m$. Assume without loss of generality that the essential variables of $f|_{\Aneq}$ are $x_1, \dotsc, x_m$. Every essential variable of $f|_{\Aneq}$ is obviously essential in $f$, so $\ess f \geq m$. Lemma \ref{existsh} implies that there is a function $h \colon A^m \to B$ such that for all $(a_1, \dotsc, a_n) \in \Aneq$ we have $f(a_1, \dotsc, a_n) = h(a_1, \dotsc, a_m)$. If a variable, say $x_n$, is inessential in $f$, then for any $(a_1, \dotsc, a_n) \in A^n$ we have $f(a_1, \dotsc, a_{n-1}, a_n) = f(a_1, \dotsc, a_{n-1}, a_{n-1}) = h(a_1, \dotsc, a_m)$. This shows that $\ess f \leq m$.
\end{proof}

\'{S}wierczkowski \cite[statement ($\beta$) in Section 2]{Sw} proves the following lemma about partitions of finite sets. Here, for any partitions $\delta$, $\delta'$ of a set $S$, we write $\delta < \delta'$ if, for each block $D$ of $\delta$, there is a block $D'$ of $\delta'$ which contains $D$.

\begin{lemma}
\label{beta}
If we have a fixed number $n \geq 3$, $S$ is a finite set, and to every partition $\delta$ of $S$ in not more than $n$ disjoint subsets corresponds a set $\varphi \delta$ of that partition so that $\delta < \delta'$ implies $\varphi \delta < \varphi \delta'$, then the intersection of all $\varphi \delta$ is nonempty.
\end{lemma}

Lemma \ref{beta} forms the base of the theorem which is usually named as \'{S}wierczkow\-s\-ki's lemma: Let $f$ be a function of arity at least $4$, such that every simple minor of $f$ is a projection. If $\delta$ is a nontrivial partition (different from equality) of the set of variables, then identifying the variables belonging to the same block of $\delta$, we obtain a simple minor $f_\delta$ of $f$. Let $\varphi \delta$ denote that block of $\delta$ which contains the variable to which $f_\delta$ projects. Then applying Lemma \ref{beta} we find that the intersection of the sets $\varphi \delta$ is not empty; thus $f_\delta$ is always the same projection, i.e., $f$ is a semiprojection.

\begin{theorem}
\label{maintheorem}
Let $f \colon A^n \to B$.
\begin{enumerate}[\rm (i)]
\item For $n \geq 2$, all variable identification minors of $f$ are constant functions if and only if $f$ is quasi-nullary.
\item For $n = 2$ or $n \geq 4$, all variable identification minors of $f$ are essentially unary if and only if $f$ is quasi-unary.
\end{enumerate}
Furthermore, in (i) and in (ii), provided that $n \geq 4$, the variable identification minors of $f$ are equivalent to the unique essentially at most unary support of $f$.
\end{theorem}
\begin{proof}
By Lemmas \ref{lem:diag} and \ref{lem:un}, it is not possible that $f$ has both constant functions and essentially unary functions as minors. Thus, in light of Remark \ref{rem1}\eqref{rem:binary}, part (ii) for $n = 2$ will follow from part (i).

(i) It is clear by definition that all minors of a quasi-nullary function are constant. For the converse implication, assume that all minors of $f$ are constant. Since $\Delta_f = \Delta_{\minor{f}{i}{j}}$ for all $i \neq j$ by Lemma \ref{lem:diag}, they must in fact be constant functions taking on the same value, say $c \in B$. Thus, $\minor{f}{i}{j}(\vect{a}) = c$ for all $\vect{a} \in A^n$ and for all $i \neq j$, so $f(\vect{a}) = c$ for all $\vect{a} \in \Aneq$. Hence $f|_{\Aneq}$ is a constant map, and so $f$ is quasi-nullary.

(ii) It is again clear by definition that all minors of a quasi-unary function are essentially unary. For the converse implication, assume that $n \geq 4$ and all minors of $f$ are essentially unary. For a nontrivial partition $\delta$ of the set of variables of $f$, denote by $f_\delta$ the simple minor of $f$ that is obtained by identifying the variables belonging to the same block. Let $\varphi \delta$ be that block of $\delta$ which contains the only essential variable of $f_\delta$. Then applying Lemma \ref{beta} we find that the intersection of the sets $\varphi \delta$ is not empty; thus $f_\delta$ always depends on the same essential variable, i.e., $f$ is quasi-unary.

The last claim follows from Lemmas \ref{lem:diag} and \ref{lem:un} and Remark \ref{rem1}\eqref{rem:unique}.
This completes the proof of Theorem \ref{maintheorem}.
\end{proof}

For $A = B$, Theorem \ref{maintheorem} restricted to semiprojections entails the well-known formulation of \'{S}wierczkowski's lemma (see, e.g., \cite{Cs,Quack}).
\begin{lemma}[\'{S}wierczkowski's lemma]
Let $f$ be an $n$-ary operation on $A$ and $n \geq 4$. Then $f$ is a semiprojection if and only if every variable identification minor of $f$ is a projection.
\end{lemma}

\section{Arity gap and a further generalization of \'{S}wierczkowski's lemma}
\label{sec:gap}

Recall that simple variable substitution induces a quasi-order on the set of $B$-valued functions on $A$, as described in Section \ref{sec:minors}. For a function $f \colon A^n \to B$ with at least two essential variables, we denote
\[
\essl f = \max_{g < f} \ess g,
\]
and we define the \emph{arity gap} of $f$ by $\gap f = \ess f - \essl f$. It is easily observed that
\[
\gap f = \min_{i \neq j} (\ess f - \ess \minor{f}{i}{j}),
\]
where $i$ and $j$ range over the set of indices of essential variables of $f$.

Since the arity gap is defined in terms of essential variables and since every $B$-valued function on $A$ is equivalent to a function whose variables are all essential, we will assume without loss of generality that the functions $f \colon A^n \to B$ whose arity gap we consider are essentially $n$-ary and $n \geq 2$.

The following upper bound for the arity gap was established by Willard \cite[Lemma 1.2]{Willard}.

\begin{theorem}
\label{Willard1.2}
Suppose $f \colon A^n \to B$ depends on all of its variables. If $n > k$, then $\gap f \leq 2$.
\end{theorem}

This theorem leaves unsettled the arity gaps of functions with a small number of essential variables, i.e., the case that $2 \leq n \leq k$.

\begin{proposition}
\label{qafltn}
Suppose $f \colon A^n \to B$, $2 \leq n \leq k$, depends on all of its variables. If $\qa f = m < n$, then $\gap f = n - m$.
\end{proposition}
\begin{proof}
Let $g$ be an essentially $m$-ary support of $f$. The variable identification minors of $f$ coincide with those of $g$, i.e., $\minor{f}{i}{j} = \minor{g}{i}{j}$ for all $i \neq j$, and hence $\ess \minor{f}{i}{j} = \ess \minor{g}{i}{j} \leq \ess g$. Since $m < n$, $g$ has an inessential variable, say $x_p$, and therefore for any $q \neq p$ we have that $\minor{f}{p}{q} = \minor{g}{p}{q} = g$. Thus, $\gap f = n - m$.
\end{proof}

In order to deal with the case that $\qa f = n$, we will adapt Willard's proof of Theorem \ref{Willard1.2} to functions of small essential arity. The idea is that it suffices to consider the restriction of $f$ to $\Aneq$ and hence the condition $n > k$ can be omitted.

\begin{theorem}
\label{GenWillard1.2}
Suppose $f \colon A^n \to B$, $n > 3$, depends on all of its variables. If $\qa f = n$, then $\gap f \leq 2$.
\end{theorem}
\begin{proof}
Let $r$ be the maximum of the essential arities of all $\minor{f}{i}{j}$, $i \neq j$. Assume on the contrary that $r < n - 2$. We shall find a contradiction.

\begin{claim}
There exist $u$, $v$ such that $\minor{f}{v}{u}$ is essentially $r$-ary and does not depend on $x_u$.
\end{claim}
Claim 1 was proved by Willard (Claim 1 in the proof of Lemma 1.2 in \cite{Willard}).

Let $u$ and $v$ be as in Claim 1, and assume without loss of generality that $u = n - 1$, $v = n$ and the essential variables of $\minor{f}{v}{u}$ are $x_1, \dotsc, x_r$. Then
\[
\minor{f}{v}{u}(x_1, \dotsc, x_n) = f(x_1, \dotsc, x_{n-2}, x_{n-1}, x_{n-1}) = h(x_1, \dotsc, x_r),
\]
where $h$ depends on all of its variables.

\begin{claim}
For all $(a_1, \dotsc, a_n) \in \Aneq$, $f(a_1, \dotsc, a_n) = h(a_1, \dotsc, a_r)$.
\end{claim}
Claim 2 was proved by Willard (Claim 2 in the proof of Lemma 1.2 in \cite{Willard}). Claim 2 implies that for all $(a_1, \dotsc, a_n) \in \Aneq$, $f|_{\Aneq}(a_1, \dotsc, a_n) = h(a_1, \dotsc, a_r)$, and hence $\qa f = \ess f|_{\Aneq} \leq r < n$, which is a contradiction.
\end{proof}

Observe that if $n > k$, then Theorem \ref{GenWillard1.2} reduces to Theorem \ref{Willard1.2}, because in this case quasi-arity means the same as essential arity.

From Proposition \ref{qafltn} and Theorem \ref{GenWillard1.2} we can now derive the following characterization of functions of arity gap at least $3$.

\begin{theorem}
\label{maintheorem2}
Suppose $f \colon A^n \to B$ depends on all of its variables. For $0 \leq m \leq n - 3$, we have that $\gap f = n - m$ if and only if $\qa f = m$.
\end{theorem}

Note that in Theorem \ref{maintheorem2} the cases $m = 0$ and $m = 1$ correspond to parts (i) and (ii) of Theorem \ref{maintheorem}, respectively. Thus, Theorem \ref{maintheorem2} can be viewed as a further generalization of \'{S}wierczkowski's lemma.

\section{The classification of pseudo-Boolean functions according to their arity gap}
\label{sec:pseudoBoolean}

In \cite{CL}, we completely classified all Boolean functions in terms of arity gap. More precisely, we have shown that for a Boolean function $f \colon \{0, 1\}^n \to \{0, 1\}$ with at least two essential variables, $\gap f = 2$ if and only if $f$ is equivalent to one of the following Boolean functions:
\begin{compactitem}
\item $x_1 + x_2 + \dots + x_m + c$ ($2 \leq m \leq n$),
\item $x_1 x_2 + x_1 + c$,
\item $x_1 x_2 + x_1 x_3 + x_2 x_3 + c$,
\item $x_1 x_2 + x_1 x_3 + x_2 x_3 + x_1 + x_2 + c$,
\end{compactitem}
where addition and multiplication are done modulo 2 and $c \in \{0, 1\}$. Otherwise $\gap f = 1$.

By Theorem \ref{Willard1.2}, the arity gap of a pseudo-Boolean function is either 1 or 2. Like in the case of Boolean functions, this fact asks for a complete classification of pseudo-Boolean functions in terms of arity gap. By making use of tools appearing in \cite{BK,Willard}, we obtain the following characterization.

\begin{theorem}
\label{gappseudoB}
For a pseudo-Boolean function $f \colon \{0,1\}^n \to B$, $n \geq 2$, which depends on all of its variables, $\gap f = 2$ if and only if $f$ satisfies one of the following conditions:
\begin{compactitem}
\item $n = 2$ and $f$ is a nonconstant function satisfying $f(0,0) = f(1,1)$,
\item $f = g \circ h$, where $g \colon \{0, 1\} \to B$ is injective and $h \colon \{0, 1\}^n \to \{0, 1\}$ is a Boolean function with $\gap h = 2$, as listed above.
\end{compactitem}
Otherwise $\gap f = 1$.
\end{theorem}

In order to prove Theorem \ref{gappseudoB}, we need to introduce some terminology and auxiliary results. Let $\mathcal{P}(A)$ denote the power set of $A$. For each positive integer $n$, define the function $\oddsupp \colon A^n \to \mathcal{P}(A)$ by
\[
\oddsupp(\vect{a}) = \{a_i : \text{$\card{\{j \in \{1, \dotsc, n\} : a_j = a_i\}}$ is odd}\}.
\]
A function $f \colon A^n \to B$ is said to be \emph{determined by $\oddsupp$} if there is a function $f^* \colon \mathcal{P}(A) \to B$ such that $f = f^* \circ \oddsupp$. In \cite{BK,Willard} it was shown that if $f \colon A^n \to B$ where $n > \max(k, 3)$ and $\gap f = 2$ then $f$ is determined by $\oddsupp$. This result leads to the following auxiliary lemma.

\begin{lemma}
\label{rangelemma}
Suppose that $f \colon A^n \to B$, where $n > \max(k, 3)$, depends on all of its variables. If the range of $f$ contains more than $2^{k-1}$ elements, then $\gap f = 1$.
\end{lemma}
\begin{proof}
For the sake of contradiction, suppose that $\gap f = 2$. Then $f$ is determined by $\oddsupp$, so $f = f^* \circ \oddsupp$ for some $f^* \colon \mathcal{P}(A) \to B$. In fact, the range of $\oddsupp$ contains only subsets of $A$ of even order or only subsets of odd order, depending on the parity of $n$. The number of subsets of $A$ of even order equals the number of subsets of odd order, and this number is $2^{k-1}$. Then the range of $f$ contains at most $2^{k-1}$ elements, contradicting the hypothesis.
\end{proof}

\begin{proof}[Proof of Theorem \ref{gappseudoB}]
It is easy to verify that if $f \colon \{0,1\}^n \to B$ satisfies any of the conditions mentioned in the statement of the theorem, then $\gap f = 2$. So we need to show that there are no other functions $f \colon \{0,1\}^n \to B$ with arity gap $2$. In fact, we only need to verify the case where the range of $f$ contains at least three elements, because otherwise $f$ is of the form $f = g \circ h$, where $g \colon \{0, 1\} \to B$ is injective and $h \colon \{0, 1\}^n \to \{0, 1\}$ is a Boolean function, and it is clear that in this case $\gap f = \gap g$.

Unary functions cannot have arity gap $2$. The case of binary functions is straightforward to verify. If $n > 3$, then Lemma \ref{rangelemma} implies that $\gap f = 1$.

If $n = 3$, we have two cases. Consider first the case that $f(0,0,0) = f(1,1,1)$. For $i \neq j$, let $g_{ij}$ be the binary function equivalent to $\minor{f}{i}{j}$. We have that $g_{ij}(0,0) = g_{ij}(1,1)$ and it is easy to see that $g_{ij}$ is either constant or essentially binary. Since $\gap f \leq 2$, the $g_{ij}$ cannot all be constant. Hence, there is an essentially binary $g_{ij}$ for some $i \neq j$, and we conclude that $\gap f = 1$.

Consider then the case that $f(0,0,0) = a \neq b = f(1,1,1)$. Let $c$ be an element in the range of $f$ distinct from both $a$ and $b$, and let $\vect{u} \in A^n$ be such that $f(\vect{u}) = c$. Then $\vect{u}$ has two coinciding coordinates, say $i$ and $j$. Let $g_{ij}$ be the binary function equivalent to $\minor{f}{i}{j}$. Then $g_{ij}$ takes on at least three distinct values, namely $a$, $b$, $c$, and it is clear that $g_{ij}$ is essentially binary, and hence $\gap f = 1$.
\end{proof}

\section{General classification of functions according to their arity gap}
\label{sec:concluding}

In the previous section, we presented an explicit classification of pseudo-Boolean functions according to their arity gap. In the general case of functions $f \colon A^n \to B$ where $\card{A} > 2$, assuming no specific algebraic structure on the underlying set $A$, such an explicit general description (in terms of representations, e.g., polynomial, DNF, etc.) seems unattainable. Falling short of explicitness, but achieving full generality, in this section we classify finite functions according to their arity gap in terms of quasi-arity and the notion of a function being determined by $\oddsupp$.

Willard's proof of Theorem 2.1 in \cite{Willard} immediately gives rise to the following generalization, where we have omitted the condition $n > \card{A}$.

\begin{theorem}
\label{thm:totsymm}
Let $f \colon A^n \to B$, $n > 3$, and suppose that $\qa f = n$ and $\gap f = 2$. Then $f|_{\Aneq}$ is totally symmetric and for all $i \neq j$, $\minor{f}{i}{j}$ depends on all of its variables except $x_i$ and $x_j$.
\end{theorem}

Berman and Kisielewicz's \cite{BK} Lemma 2.7 can similarly be generalized, removing the condition $n > \card{A}$, which together with Theorem \ref{thm:totsymm} and Proposition \ref{qafltn} yields the following result providing a necessary and sufficient condition for a function to have arity gap $2$. For $f \colon A^n \to B$, we say that \emph{$f|_{\Aneq}$ is determined by $\oddsupp$} if $f|_{\Aneq} = f^* \circ \oddsupp$, where $f^* \colon \mathcal{P}'(A) \to B$ is a nonconstant function and  $\oddsupp \colon \Aneq \to \mathcal{P}'(A)$ is defined as in Section \ref{sec:pseudoBoolean}, but here $\mathcal{P}'(A)$ denotes the set of odd or even---depending on the parity of $n$---subsets of $A$ of order at most $n - 2$.

\begin{theorem}
\label{thm:strongly}
Suppose $f \colon A^n \to B$, $n > 3$, depends on all of its variables. Then $\gap f = 2$ if and only if $\qa f = n - 2$ or $\qa f = n$ and $f|_{\Aneq}$ is determined by $\oddsupp$.
\end{theorem}

We now have everything ready to present a complete classification of all functions $f \colon A^n \to B$ according to their arity gap.

\begin{theorem}
\label{thm:last}
Suppose that $f \colon A^n \to B$, $n \geq 2$, depends on all of its variables.
\begin{compactenum}[\rm (i)]
\item For $0 \leq m \leq n - 3$, $\gap f = n - m$ if and only if $\qa f = m$.
\item For $n \neq 3$, $\gap f = 2$ if and only if $\qa f = n - 2$ or $\qa f = n$ and $f|_{\Aneq}$ is determined by $\oddsupp$.
\item For $n = 3$, $\gap f = 2$ if and only if there is a nonconstant unary function $h \colon A \to B$ and $i_1, i_2, i_3 \in \{0, 1\}$ such that
\begin{align*}
& f(x_1, x_0, x_0) = h(x_{i_1}), \\
& f(x_0, x_1, x_0) = h(x_{i_2}), \\
& f(x_0, x_0, x_1) = h(x_{i_3}).
\end{align*}
\item Otherwise $\gap f = 1$.
\end{compactenum}
\end{theorem}

\begin{proof}
Statements (i) and (ii) follow from Theorems \ref{maintheorem2} and \ref{thm:strongly}, respectively.

For (ii), the sufficiency is obvious. In order to prove the necessity, assume that $f \colon A^3 \to B$ has arity gap $2$. Then all variable identification minors of $f$ have essential arity at most $1$ and at least one of them has essential arity $1$. By Lemma \ref{lem:diag}, all of them have essential arity $1$ and are actually of the form $\Delta_f(x_{i_j})$ for some $i_1, i_2, i_3 \in \{1, 2, 3\}$.

Statement (iv) is clear, because (i)--(iii) exhaust all other possibilities.
\end{proof}

Let $f$ be a ternary function with arity gap $2$, and let $h$, $i_1$, $i_2$, $i_3$ be as in Theorem \ref{thm:last}(iii). If $(i_1, i_2, i_3) = (1, 0, 0)$, $(0, 1, 0)$ or $(0, 0, 1)$, then $\qa f = 1$; in the other five cases $\qa f = 3$. If $(i_1, i_2, i_3) = (1, 1, 1)$, then $f|_{\Aneq}$ is determined by $\oddsupp$, otherwise not. Thus $(i_1, i_2, i_3) = (1, 1, 0)$, $(1, 0, 1)$, $(0, 1, 1)$ and $(0, 0, 0)$ give counterexamples to show that statement (ii) of Theorem \ref{thm:last} does not hold when $n = 3$ (and these are the only counterexamples).

Note that in the case $A = B$ and $h = \id_A$ we have that $(i_1, i_2, i_3) = (0, 0, 0)$ if and only if $f$ is a majority operation; $(i_1, i_2, i_3) = (1, 1, 1)$ if and only if $f$ is a minority operation; $(i_1, i_2, i_3) = (1, 0, 0)$, $(0, 1, 0)$ or $(0, 0, 1)$ if and only if $f$ is a semiprojection; and $(i_1, i_2, i_3) = (1, 1, 0)$, $(1, 0, 1)$ or $(0, 1, 1)$ if and only if $f$ is a so-called $\frac{2}{3}$-minority operation.

\section*{Acknowledgements}

We would like to thank Jorge Almeida, Michael Pinsker, Maurice Pouzet, and Ross Willard for useful discussions on the topic.

Part of this work was carried out while both authors were visiting Tampere and while the second author was visiting the University of Luxembourg. We would like to thank the Department of Mathematics, Statistics and Philosophy of the University of Tampere and the Department of Mathematics of the University of Luxembourg for providing working facilities.

\end{document}